\documentclass[11pt]{article}
\usepackage{amsmath}
\usepackage{amssymb}

\newtheorem{thm}{Theorem}[section]
\newtheorem{lemma}[thm]{Lemma}
\newtheorem{prop}[thm]{Proposition}
\newtheorem{rem}[thm]{Remark}
\newtheorem{cor}[thm]{Corollary}
\newtheorem{defn}[thm]{Definition}

\begin{document}

\title{Completely
    monotonic functions of positive order and asymptotic expansions of the
    logarithm of Barnes double gamma
    function and Euler's gamma function}
\author{Stamatis Koumandos and Henrik L. Pedersen}
\date{\today}
\maketitle


\begin{abstract}{We introduce completely monotonic functions of order
    $r>0$ and show that the remainders in asymptotic expansions of the
    logarithm of Barnes double gamma function and Euler's gamma
    function give rise to completely monotonic functions of any
    positive integer order.}
\end{abstract}

\section{Completely monotonic functions of positive order}
In this paper it is established that the remainders in asymptotic
expansions of the logarithm of Barnes double
gamma function and Euler's gamma function are (up to a sign) 
completely monotonic functions of order comparable with the decay of
the remainder in the expansion. 

We first recall some definitions and give some preliminary results.

A function $f\,:\,(0,\,\infty) \rightarrow \mathbb{R}$ is called
completely monotonic if $f$ has derivatives of all orders and
satisfies
\begin{equation}
\label{eq:cm}
(-1)^n\,f^{(n)}(x) \geq 0,\;\;\textrm{for
all}\;\;x>0\;\;\textrm{and}\;\;n=0,1,2,\ldots
\end{equation}
J.\ Dubourdieu \cite{du} proved that if a non constant function $f$
is completely monotonic then strict inequality holds in \eqref{eq:cm}. See
also \cite{van} for a simpler proof of this result. A
characterization of completely monotonic functions is given by
Bernstein's theorem, see \cite[p.\ 161]{w}, which states that $f$ is
completely monotonic if and only if
$$
f(x)=\int _{0}^{\infty} e^{-xt}\,d\mu(t),
$$
where $\mu$ is a nonnegative measure on $[0,\,\infty)$ such that
the integral converges for all $x>0$.

Here we are interested in the class of strongly completely monotonic
functions, introduced in \cite{tww}. A function $f\,:\,(0,\,\infty)
\rightarrow \mathbb{R}$ is called strongly completely monotonic if
it has derivatives of all orders and
$(-1)^{n}\,x^{n+1}\,f^{(n)}(x)$ is nonnegative and decreasing on
$(0,\,\infty)$ for all $n=0,1,2,\ldots$. (It is clear that being strongly
completely monotonic is {\em stronger} that being completely monotonic.)
These functions are connected to the important question of superadditivity (cf.\ \cite{tww}).

The following proposition contains a simple characterization of
strongly completely monotonic functions and its proof is obtained by a
direct application of the definitions given above. 
\begin{prop}
\label{prop:simple_char}
A function $f(x)$ is strongly completely monotonic if and only
  if the function $xf(x)$ is completely monotonic.
\end{prop}
In \cite{tww} the authors gave another characterization of
strongly completely monotonic functions.
\begin{prop}
\label{prop:char}
The function $f(x)$ is strongly completely monotonic if and only if
$$
f(x)=\int_{0}^{\infty} e^{-xt}\,p(t)\,dt,
$$
where $p(t)$ is nonnegative and increasing and the integral
converges for all $x>0$.
\end{prop}
We notice that any right-continuous function $p$ appearing in Proposition \ref{prop:char} can be written as
$$
p(t)=\mu([0,t]),
$$
where $\mu$ is a Radon measure on $[0,\infty)$.

An extension of Proposition \ref{prop:simple_char} and \ref{prop:char} is the following.
\begin{thm}
\label{thm:r_char}
Let $r$ be an integer $\geq 1$. The function
$x^r\,f(x)$ is completely monotonic if and only if
\begin{equation*}
f(x)=\int_{0}^{\infty} e^{-xt}\,p(t)\,dt\,,
\end{equation*}
where the integral converges for all $x>0$ and where $p$ is $r-1$ times differentiable on $[0,\infty)$ with 
$p^{(r-1)}(t)=\mu ([0,t])$ for some Radon measure $\mu$ and
$p^{(k)}(0)=0$ for $0\leq k \leq r-2$. 
\end{thm}
\begin{rem}
The conditions on $p$ ensure that $p, p', \ldots ,p^{(r-1)}$ are all nonnegative.
\end{rem}
{\em Proof}.
Suppose that $x^rf(x)$ is completely monotonic. Then  
$$
x^rf(x)=\int_0^{\infty}e^{-xt}\, d\mu (t)={\cal L} (\mu)(x),
$$
for some Radon measure $\mu$ (and where $\cal L$ denotes the Laplace transform). Furthermore, it is a fact that 
$$
x^{-r}=\frac{1}{(r-1)!}\int_0^{\infty}t^{r-1}e^{-xt}\, dt,
$$
and this yields $f(x)={\cal L}(p)(x)$ where $p$ is the convolution of
these two measures on the half line,
$$
p(t)=\left(\frac{s^{r-1}}{(r-1)!}ds\ast \mu\right)
(t)=\frac{1}{(r-1)!}\int_0^{t}(t-s)^{r-1}\, d \mu (s). 
$$
From this formula it is easy to check that $p$ has the asserted properties.

Conversely, if 
$$
p^{(r-1)}(t)=\mu ([0,t])
$$
and all derivatives of $p$ at $t=0$ up to order $r-2$ are zero then we find by
integration ($r-1$ times)
$$
p(t)=\frac{1}{(r-1)!}\int_0^{t}(t-s)^{r-1}\, d \mu (s).  
$$
This gives
\begin{eqnarray*}
x^rf(x)&=&x^r
{\cal L}(p)(x)\\
&=&x^r{\cal L}\left(\frac{s^{r-1}}{(r-1)!}ds\ast \mu\right)(x)\\
&=&x^r{\cal L}\left(\frac{s^{r-1}}{(r-1)!}ds\right)(x){\cal L}\left(\mu\right)(x)\, =\, {\cal L}\left(\mu\right)(x),
\end{eqnarray*}
and therefore $x^rf(x)$ is completely monotonic.\hfill $\square$

In the light of these results we formulate the following definition.
\begin{defn}
Let $r\geq 0$. A function $f$ defined on $(0,\infty)$ is
said to be completely monotonic of order $r$ if $x^rf(x)$ is
completely monotonic.
\end{defn}
According to this definition,  completely monotonic functions of order
0 are the classical completely monotonic functions, order 1 are the
strongly completely monotonic functions  and so on. 

Theorem \ref{thm:r_char} is the basis of the proofs of Theorems
\ref{thm:euler_strong}, \ref{thm:G} and \ref{thm:main1}.

\begin{rem}
Following \cite[Section 2.9]{aar} a fractional integral
$I_{\alpha}(\mu )(t)$ (for $\alpha>0$) of a measure  on 
$[0,\infty)$ is defined by 
$$
I_{\alpha}(\mu )(t)=\frac{1}{\Gamma(\alpha)}\int_0^{t}(t-s)^{\alpha-1}\, d \mu (s). 
$$
The proof of Theorem \ref{thm:r_char} actually shows that
$f$ is completely monotonic of order $\alpha >0$ if and only if $f$
is the Laplace transform of a fractional integral of a positive Radon
measure on $[0,\infty)$, that is,
$$
f(x)=\frac{1}{\Gamma(\alpha)}\int_0^{\infty}e^{-xt}\int_0^{t}(t-s)^{\alpha-1}\,
d \mu (s)\, dt.
$$
\end{rem}
\bigskip

\section{Euler's gamma function}
Let us consider the asymptotic expansion of the logarithm of Euler's
gamma function
\begin{eqnarray}
&&\log\,\Gamma(x)=\Big(x-\frac{1}{2}\Big)\log
x-x+\frac{1}{2}\log(2\pi) +\nonumber\\[+8pt] &&
+\sum_{k=1}^{n}\frac{B_{2k}}{(2k-1)\,2k}\,\frac{1}{x^{2k-1}}+(-1)^n
\,R_{n}(x),\label{eq:euler-remainder}
\end{eqnarray}
where $B_{2k}$ are the Bernoulli numbers.
It has been shown in \cite{al} that for all $n=0,1,2,\ldots$,
the remainder $R_{n}(x)$ is completely monotonic on $(0,\,\infty)$.
We strengthen this result in Theorem \ref{thm:euler_strong} below.

For $t>0$ we write
\begin{equation}
\label{eq:def_V_n}
\frac{t}{e^t-1}=1-\frac{t}{2}+\sum_{k=1}^{n}\frac{B_{2k}}{(2k)!}\,t^{2k}+(-1)^n\,t^{2n+2}\,V_{n}(t), 
\end{equation}
where the remainder term $V_{n}(t)$
is given as (see \cite[p. 64]{te})
\begin{equation}
\label{eq:V_n}
V_{n}(t)=\sum_{k=1}^{\infty}\frac{2}{(t^2+4\pi^2\,k^2)(2\pi\,
k)^{2n}}, \quad n\geq 0.
\end{equation}

\begin{thm}
\label{thm:euler_strong}
The function $R_{n}(x)$ defined in \eqref{eq:euler-remainder} is
completely monotonic of order $k$ for $n\geq k$ for any $k\geq 0$. Indeed,
$$
R_{n}(x)=\int_{0}^{\infty}e^{-xt}\,r_{n}(t)\,dt,
$$
where $r_{n}(t)=t^{2n}V_n(t)$ satisfies $r_{n}^{(l)}(0)=0$ for $0\leq
l\leq k-1$ and $r_{n}^{(k)}(t)>0$ for $t>0$.
\end{thm}
A similar result regarding complete monotonicity of the remainder of  an asymptotic expansion of a ratio of gamma functions of the form
$\Gamma(x+a)/\Gamma(x+b)$ in terms of powers of $1/(x+c)$,  is given
in Frenzen's  paper \cite{f}.
\begin{rem}
In \cite{k} it is proved that $V_{n}(t)$ takes
also the form
$$
V_{n}(t)=\frac{1}{(2n+1)!}\,\frac{1}{e^t-1}\,\int_{0}^{1}e^{tu}\,(-1)^n\,B_{2n+1}(u)\,du,\quad
n\geq 0,
$$
where $B_{2n+1}(u)$ are the Bernoulli polynomials. We observe that this formula extends
to $t=0$ because of the well-known property
$$
\int_0^1B_{2n+1}(u)\, du=0
$$
for $n\geq 0$.
\end{rem}
In \cite{k2} it is proved that $V_n$ is positive, decreasing and
satisfies $(t^2\,V_{n}(t))^{\prime}>0$. 
We extend the last property in the
following Lemma \ref{lemma:V_n_ny}.  
\begin{lemma} 
\label{lemma:V_n_ny}
The function $V_n$ has the following properties:
\begin{enumerate} 
\item[(i)]
$
(t^{2j}V_n(t))^{(l)}>0
$
for all $n\geq 0$ and all $l\leq j$,
\item[(ii)]
$
(t^{2j-1}V_n(t))^{(l)}>0
$
for all $n\geq 0$ and all $l\leq j-1$.
\end{enumerate}
\end{lemma}
{\it Proof.} We have 
$$
t^{2j}V_n(t)=2\sum_{k=1}^{\infty}\frac{t^{2j}}{t^2+(2\pi k)^2}\frac{1}{(2\pi
  k)^{2n}}=2\sum_{k=1}^{\infty}\frac{1}{(2\pi k)^{2n-2j+2}}s_j(t/(2\pi k)),
$$
where 
$$
s_j(x)=\frac{x^{2j}}{1+x^2}.
$$
Here we notice that  $s_j^{(l)}(x)$ tends to zero as $x^{2j-l}$ for $x$
tending to 0. We get, by differentiation,
$$
(t^{2j}V_n(t))^{(l)}=2\sum_{k=1}^{\infty}\frac{1}{(2\pi k)^{2n-2j+2}}s_j^{(l)}(t/(2\pi k))\frac{1}{(2\pi k)^l},
$$
where the series converges uniformly due to the behaviour of
$s_j^{(l)}(t/(2\pi k))$ for large $k$.  We have $s_j(x)=\xi_j(x^2)$,
where
$$
\xi_j(x)=\frac{x^j}{1+x}
$$
and from Lemma
\ref{lemma:xi} below it follows that $s_j^{(l)}(x)$ are all
positive. This proves the first assertion.  The
second assertion is proved in the same way. (Of course the first
assertion follows from the second when $l\leq j-1$.) \hfill $\square$
\begin{lemma}
\label{lemma:xi}
Let 
$$
\xi_n(x)=\frac{x^n}{1+x}.
$$
Then $\xi_n^{(k)}(x)>0$ for $k\leq n$ and $x>
0$. 
\end{lemma}
{\it Proof.}
First of all we notice that $\xi_n^{(k)}(0)=0$ for $k=0,\ldots,
n-1$. To consider the $n$'th derivative, $\xi_n(x)$ is rewritten 
as follows:
\begin{eqnarray*}
\xi_n(x)&=&\frac{x^{n}}{1+x}\\
&=&\frac{1}{x+1}\sum_{k=0}^n\binom{n}{k}(x+1)^k(-1)^{n-k}\\
&=&\frac{(-1)^n}{x+1}+\sum_{k=1}^n\binom{n}{k}(x+1)^{k-1}(-1)^{n-k}
\end{eqnarray*}
so that 
$$
\xi_n^{(n)}(x)=(-1)^{n}(x+1)^{-(n+1)}(-1)^nn!=\frac{n!}{(x+1)^{n+1}}>0.
$$
Using the relation 
$$
\xi_n^{(k)}(x)=\int_0^{x}\xi_n^{(k+1)}(t)\, dt
$$
recursively for $k$ from $n-1$ to 0, we find that the derivatives of
$\xi_n$ of order not exceeding $n$ are positive.\hfill $\square$  

In the next proposition we gather some additional properties of the
function $V_n$ that are of interest in their own right.  
\begin{prop} 
\label{prop:V_n_additional}
The function $V_n$ has the following additional properties.
\begin{enumerate} 
\item[(i)]
  $t^{2n}V_n(t)=(-1)^nV_0(t)+\sum_{k=1}^{n}(-1)^{n-k}V_{n-k}(0)t^{2k-2}$,\;\; $n\geq 0$.
\item[(ii)] $V_{n}(0)=(-1)^{n}\,\frac{B_{2n+2}}{(2n+2)!}$,\,\, $n\geq 0$.
\item[(iii)] $V_{n}^{\prime\prime}(0)=-2\,V_{n+1}(0)$,\,\, $n\geq 0$. 
\item[(iv)] $V_{n}^{\prime\prime}(t)-V_{n}^{\prime\prime}(0)>0$,$\;\;
\forall\, t>0,\;\;n\geq 0$.
\end{enumerate}
\end{prop}
\begin{rem}
The function $V_0$ appearing in this proposition satisfies by definition $t/(e^t-1)= 1-t/2+t^2
V_{0}(t)$, whence $V_{0}(t)=1/t^2((t/2)\coth(t/2)-1)$.
\end{rem}
{\it Proof.} The first assertion is obtained by repeated application
of the recursive relation 
$$
t^2\,V_{n}(t)=V_{n-1}(0)-V_{n-1}(t),
$$
see \cite{k2}. The remaining
assertions are easily obtained by using definition
\eqref{eq:def_V_n} and relation \eqref{eq:V_n}.\hfill $\square$
 
\noindent
 \emph{Proof of Theorem \ref{thm:euler_strong}}. Using Binet's formula
\begin{eqnarray*}
&&\log\,\Gamma(x)=\Big(x-\frac{1}{2}\Big)\log
x-x+\frac{1}{2}\log(2\pi) +\\[+8pt]
&&+\int_{0}^{\infty}\Big(\frac{t}{2}-1+\frac{t}{e^{t}-1}\Big)\,\frac{e^{-xt}}{t^2}\,dt,\;\;x>0
\end{eqnarray*}
and formula \eqref{eq:def_V_n} we see that
$$
r_{n}(t)=t^{2n}\,V_{n}(t),\;\;n=0,1,2,\ldots
$$
It is clear that $r_{n}^{(k)}(0)=0$ for $k\leq n-1$. Now use $(i)$ of Lemma \ref{lemma:V_n_ny} to prove that
$r_{n}^{(n)}(t)>0$ for $t>0$. The result then follows from 
Theorem \ref{thm:r_char}. \hfill $\Box$

\begin{cor}
For any $n\geq 1$ the function
\begin{eqnarray*}
F_n(x)&=&(-1)^n\;\left[x^n\,\log\,\Gamma(x)-x^n\,\left(x-\frac{1}{2}\right)\log
x+x^{n+1}-\frac{x^n}{2}\,\log(2\pi)\right.\\
&&\left. -\sum_{k=1}^{n}\frac{B_{2k}}{(2k-1)\,2k}\,x^{n-2k+1}\right]
\end{eqnarray*}
is completely monotonic on $(0, \infty)$.
\end{cor}
{\it Proof.} It follows by a combination of Theorem \ref{thm:euler_strong}
with Theorem \ref{thm:r_char} and \eqref{eq:euler-remainder}. \hfill $\Box$

\section{Barnes $G$-function}
We note that a result similar to Theorem \ref{thm:euler_strong} holds for the remainder
in an asymptotic expansion (due to C.~Ferreira and J.~L.~L\'{o}pez,
\cite[Theorem 1]{f1}) of the logarithm of Barnes $G$-function. This
function is defined as an infinite product and satisfies 
$G(1)$ =1 and
$G(z+1)=\Gamma(z)\,G(z)$. See also \cite{ped} for details and
additional considerations. The remainder in this expansion takes the form
\begin{equation*}
P_{n}(x)=(-1)^n\,\int_{0}^{\infty}e^{-xt}\,t^{2n-1}\,V_{n}(t)\,dt,
\end{equation*}
where $V_{n}(t)$ is as above.
\begin{thm}
\label{thm:G}
The remainder $(-1)^{n} P_{n}(x)$ is completely monotonic of order
$k$ on $(0,\,\infty)$ for $n\geq k+1$.
\end{thm} 
{\it Proof.} Let $\lambda_n(t)=t^{2n-1}\,V_{n}(t)$. Clearly,
$\lambda_n^{(k)}(0)=0$ for $k\leq n-1$. It follows from Lemma
\ref{lemma:V_n_ny} that $\lambda_n^{(n-1)}(t)>0$ for $t>0$, and the
result follows from Theorem \ref{thm:r_char}.\hfill
$\square$ 

\bigskip

\section{Barnes double gamma function}
This section is devoted to the investigation of the remainders in an
asymptotic expansion due to Ruijsenaars of the logarithm of Barnes
double gamma function. The expansion is given in terms of generalized
Bernoulli polynomials $B_{k}^{(2)}(x)$, see \cite[p.\,615]{aar} or
\cite[p.\,4]{te}. Our investigation is based on Ruijesnaars' results and
therefore we have found it natural to his terminology
$B_{2,k}(x)=B_{k}^{(2) }(x)$. We shall futhermore call these polynomials
the double Bernoulli polynomials. (Ruijsenaars calls $B_{N,k}(x)$
multiple Bernoulli polynomials.)

The double Bernoulli polynomials $B_{2,k}(x)$ are defined by
$$
\frac{t^2e^{xt}}{(e^{t}-1)^2}= \sum_{k=0}^{\infty}B_{2,k}(x)\frac{t^k}{k!}
$$
and the double Bernoulli numbers $B_{2,k}$ by
$B_{2,k}=B_{2,k}(0)$.

The asymptotic expansion of the logarithm of Barnes double
gamma function with both parameters equal to 1, $\log \Gamma_2(w)=\log
\Gamma_2(w|1,1)$, is given as follows:

\begin{eqnarray*}
\log \Gamma_2(w)&=& -\frac{B_{2,2}(w)}{2}\log w +
\frac{3}{4}B_{2,0}w^2+B_{2,1}w\\
&& +\sum_{k=3}^{M}\frac{(-1)^k}{k!}(k-3)!B_{2,k}w^{2-k} + R_{2,M}(w),
\end{eqnarray*}
where the remainder $R_{2,M}$ has the representation
$$
R_{2,M}(w)= \int_{0}^{\infty}\frac{e^{-wt}}{t^3}\left(
\frac{t^2}{(1-e^{-t})^2}-
    \sum_{k=0}^{M}\frac{(-1)^k}{k!}B_{2,k}t^k\right) \, dt.
$$
Here, $\Re w>0$ and $M\geq 2$. See \cite[(3.13) and (3.14)]{rui}.




In \cite{k2} and \cite{ped2} it was shown independently that
$(-1)^{n-1}R_{2,2n}(x)$ is a completely monotonic function. Below it is
verified that it is indeed a completely monotonic function of order
$k$ for
$n\geq k+1$.

We briefly indicate the two different proofs of complete
monotonicity. By Bernstein's theorem it amounts to showing the
positivity of $U_n(t)$ for $t>0$ and any $n\geq 1$, where
$$
U_n(t)=(-1)^{n-1}\left( \frac{t^2}{(1-e^{-t})^2}-
    \sum_{k=0}^{2n}\frac{(-1)^k}{k!}B_{2,k}t^k\right),
$$
since
\begin{equation}
\label{eq:kernel}
(-1)^{n-1}R_{2,2n}(x)= \int_{0}^{\infty}e^{-xt}\frac{U_n(t)}{t^3}\, dt.
\end{equation}
In \cite{k2} it was shown that
\begin{equation}
\label{eq:U_n}
U_n(t)=t^{2n+1}V_{n-1}(t)+t^2\left(t^{2n+1}V_n(t)\right)',
\end{equation}
where $V_{n}(t)$ is defined in \eqref{eq:def_V_n} and the proof is obtained by showing that
\newline $(t^{2n+1}V_n(t))'>0$.
In \cite{ped2} the proof is based on a contour integration argument
and the following representation of $p_n(t)=U_n(t)/t^3$
in \eqref{eq:kernel} is found
\begin{equation}
\label{eq:p_n}
\begin{array}{rl}
p_n(t)=\displaystyle{t^{2n-2}\sum_{k=1}^{\infty}(2\pi k)^{1-2n}} &\left(  \displaystyle{\frac{4\pi k}{t^2+(2\pi k)^2}+\frac{8\pi kt}{(t^2+(2\pi k)^2)^2}+}\right. \\
& \left.  \displaystyle{\frac{(2n-1)}{2\pi k}\frac{2t}{t^2+(2\pi
    k)^2}}\right).
\end{array}
\end{equation}
This clearly shows the positivity of $U_n$.

The main result is formulated in the theorem below.
\begin{thm}
\label{thm:main1} 
The remainder $(-1)^{n-1}\,R_{2,2n}(x)$ is completely monotonic of
order $k$ on $(0,\,\infty)$ for $n\geq k+1$.
\end{thm}
{\it Proof.} It follows from
\eqref{eq:U_n} that
\begin{eqnarray*}
 p_{n}(t)&=&\frac{U_{n}(t)}{t^3}=t^{2n-2}V_{n-1}(t)+(2n+1)t^{2n-1}V_{n}(t)+t^{2n}V^{\prime}_{n}(t)\\[+3pt]
 &=&t^{2n-2}V_{n-1}(t)+t^{2n-1}V_{n}(t)+(t^{2n}V_{n}(t))^{\prime}\\[+3pt]
 &=&r_{n-1}(t)+\lambda_{n}(t)+r^{\prime}_{n}(t),
\end{eqnarray*}
where $r_{n}(t)=t^{2n}\,V_{n}(t)$ and $\lambda_{n}(t)=t^{2n-1}\,V_{n}(t)$.

Clearly $p^{(k)}_{n}(0)=0$ for $k\leq n-2$. Since $r^{(k)}_{n}(t)>0$ for $0\leq k\leq n$ and
$\lambda^{(k)}_{n}(t)>0$ for $0\leq k\leq n-1$, $t>0$,  it follows
from the above that
 $p^{(k)}_{n}(t)>0$ for $0\leq k\leq n-1$ and for
all $t>0$ and by Theorem \ref{thm:r_char} this completes the proof of the Theorem. \hfill $\square$
\begin{rem}
Theorem \ref{thm:main1} can also be obtained directly from the
representation \eqref{eq:p_n}. Indeed, since
$$
t^{2n-2}\sum_{k=1}^{\infty}(2\pi k)^{1-2n}\frac{4\pi k}{t^2+(2\pi k)^2}=
\sum_{k=1}^{\infty}\frac{2}{(2\pi k)^2}\frac{\left( t/(2\pi k)\right)^{2n-2}}{1+\left(t/(2\pi k)\right)^2},
$$
$$
t^{2n-2}\sum_{k=1}^{\infty}(2\pi k)^{1-2n} \frac{8\pi kt}{(t^2+(2\pi k)^2)^2}=
\sum_{k=1}^{\infty}\frac{4t}{(2\pi k)^4}\frac{\left( t/(2\pi k)\right)^{2n-2}}{(1+\left(t/(2\pi k)\right)^2)^2}
$$
and
$$
t^{2n-2}\sum_{k=1}^{\infty}(2\pi k)^{1-2n} \frac{(2n-1)}{2\pi k}\frac{2t}{t^2+(2\pi k)^2}=
\sum_{k=1}^{\infty}\frac{2(2n-1)t}{(2\pi k)^4}\frac{\left( t/(2\pi k)\right)^{2n-2}}{1+\left(t/(2\pi k)\right)^2}
$$
we have from \eqref{eq:p_n}
$$
p_n(t)=\sum_{k=1}^{\infty}\frac{2}{(2\pi k)^2}s_{n-1}\left( t/(2\pi k)\right)+
\sum_{k=1}^{\infty}\frac{4t}{(2\pi k)^4}g_{n-1}\left( t/(2\pi k)\right),
$$
where $g_n$ is defined by 
$$
g_n(x)=(n+1/2)\frac{x^{2n}}{1+x^2}+\frac{x^{2n}}{(1+x^2)^2}.
$$
We put 
$$
h_n(x)=n\frac{x^{n}}{1+x}+\frac{x^{n}}{(1+x)^2}
$$
and have in this way $g_n(x)=s_{n}(x)/2+h_n(x^2)$. The positivity of
the $n$'th
derivative of $g_n$ clearly follows from the positivity of the
derivatives $h_n^{(k)}$ for $k\leq n$. To investigate the derivatives
of $h_n$ we rewrite it as follows.
\begin{eqnarray*}
h_n(x)&=&\frac{nx^n}{x+1}+\frac{x^n}{(x+1)^2}\\
&=&\frac{n}{x+1}\left\{(-1)^n+\sum_{k=1}^n\binom{n}{k}(x+1)^k(-1)^{n-k}\right\}\\
&&+\frac{1}{(x+1)^2}\Bigg\{(-1)^n+n(x+1)(-1)^{n-1}+\Big.\\
&&\left.+\sum_{k=2}^n\binom{n}{k}(x+1)^k(-1)^{n-k}\right\}\\
&=&\frac{(-1)^n}{(x+1)^2}+l_n(x),
\end{eqnarray*}
where $l_n$ is a polynomial of degree $n-1$. Therefore
$$
h_n^{(n)}(x)= \frac{(n+1)!}{(1+x)^{n+2}}>0.
$$
Furthermore, $h_n^{(k)}(0)=0$ for $k\leq n-1$, whence 
$$
h_n^{(k-1)}(x)=\int_0^xh_n^{(k)}(t)dt>0
$$
for $k=1,\ldots, n$.
This completes a different proof of Theorem \ref{thm:main1}.
\end{rem}

\noindent
{\it Acknowledgment.} The authors thank Christian Berg for his
comments in particular regarding Theorem \ref{thm:r_char}.

\vspace{0.8cm}

\noindent
Stamatis Koumandos\\
Department of Mathematics and Statistics\\[-2pt]
The University of Cyprus \\[-2pt]
P. O. Box 20537\\[-2pt]
1678 Nicosia\\[-2pt]
 CYPRUS\\[-1pt]
{\em email}:\hspace{2mm}{\tt skoumand@ucy.ac.cy}

\vspace{0.5cm}

\noindent
Henrik Laurberg Pedersen\\
Department of Basic Sciences and Environment\\[-2pt]
Mathematics and Computer Science\\[-2pt]
Faculty of Life Sciences\\[-2pt]
University of Copenhagen \\[-2pt]
40, Thorvaldsensvej\\[-2pt]
DK-1871 Frederiksberg C\\[-2pt]
DENMARK\\[-1pt]
{\em email}:\hspace{2mm}{\tt henrikp@dina.kvl.dk}

\end{document}